\documentclass[a4paper,reqno,10pt]{amsart}

\usepackage{amsmath}
\usepackage{amsthm}
\usepackage{amssymb}
\usepackage{amscd} 

\usepackage{mathabx}

\usepackage{bbm} 
\usepackage{mathrsfs} 

\usepackage[applemac]{inputenc} 

\usepackage{verbatim} 

\swapnumbers 

\newtheoremstyle{exercise} 
  {3pt} 
  {3pt} 
  {\scriptsize\rmfamily} 
  {
\parindent} 
  {\rmfamily\scshape} 
  {.} 
  {.5em} 
  {} 

\newtheoremstyle{newplain}
  {5pt}
  {5pt}
  {\itshape}
  {}
  {\rmfamily\scshape}
  {. ---}
  {.5em}
  {}

\newtheoremstyle{newremark}
  {5pt}
  {5pt}
  {\rmfamily}
  {}
  {\rmfamily\scshape}
  {. ---}
  {.5em}
  {}

\theoremstyle{newplain}
\newtheorem*{Theorem*}{Theorem} 

\theoremstyle{newplain}
\newtheorem{Theorem}{Theorem}
\newtheorem{Lemma}[Theorem]{Lemma}

\newtheorem{Proposition}[Theorem]{Proposition}
\newtheorem{Conjecture}[Theorem]{Conjecture}
\newtheorem{Definition}[Theorem]{Definition}

\theoremstyle{newremark}
\newtheorem{Remark}[Theorem]{Remark}
\newtheorem{Question}[Theorem]{Question}
\newtheorem{Example}[Theorem]{Example}
\newtheorem{Claim}[Theorem]{Claim}

\theoremstyle{exercise}

\numberwithin{Theorem}{section}
\numberwithin{Exercise}{subsection}

\newcommand{\R}{\mathbb{R}} 

\newcommand{\calA}{\mathscr{A}}

\newcommand{\calD}{\mathscr{D}}
\newcommand{\calE}{\mathscr{E}}

\newcommand{\calH}{\mathscr{H}}

\newcommand{\calL}{\mathscr{L}}
\newcommand{\calM}{\mathscr{M}}

\newcommand{\calP}{\mathscr{P}}

\newcommand{\calW}{\mathscr{W}}

\def\XXint#1#2#3{{%
\setbox0=\hbox{$#1{#2#3}{\int}$}
\vcenter{\hbox{$#2#3$}}\kern-.5\wd0}}

\renewcommand{\geq}{\geqslant}
\renewcommand{\subset}{\subseteq}

\begin{document}
\title[Hydrodynamic flow]{Recent developments of analysis for hydrodynamic flow of nematic liquid
crystals}

\begin{abstract}
The study of hydrodynamics of liquid crystal leads to many fascinating mathematical problems,
which has prompted various interesting works recently.
This article reviews the static Oseen-Frank theory and surveys some recent progress on the existence,
regularity, uniqueness, and large time asymptotic of the hydrodynamic flow of nematic liquid crystals.
We will also propose a few interesting questions for future investigations.
\end{abstract}

\author{Fanghua Lin}
\address{Courant Institute of Mathematical Sciences,
New York University, NY 10012, USA
and NYU-ECNU Institute of Mathematical Sciences, at NYU Shanghai, 3663, North Zhongshan Rd., Shanghai,
P. R. China 200062}
\email{linf@cims.nyu.edu}
\author{Changyou Wang}
\address{
Department of Mathematics,
University of Kentucky,
Lexington, KY 40506,
USA, and Department of Mathematics, Purdue University, 150 
N. University Street, West Lafayette, IN 47907, USA}
\email{wang2482@purdue.edu}
\maketitle

\setcounter{section}{0} \setcounter{equation}{0}
\section{Static Theory}

Liquid crystal is an intermediate phase between the crystalline solid state
and the isotropic fluid state. It possess none or partial positional order
but displays an orientational order at the mean time.
There are roughly three types of liquid crystals commonly referred in the
literature: nematic, cholesteric, and smectic. The nematic liquid crystals are composed of
rod-like molecules with the long axes of neighboring molecules approximately aligned to
one another. The simplest continuum model to
study the equilibrium phenomena for nematic liquid crystals is the Oseen-Frank
theory, proposed by Oseen \cite{Oseen}  in 1933 and Frank \cite{Frank} in 1958.
A unifying program describing general liquid crystal materials is the Landau-de Gennes theory (\cite {LDG1}
and \cite{LDG2}), which involves the orientational order parameter tensors.

In the absence of surface energy terms and applied fields, the Oseen-Frank theory
seeks unit-vector fields $d:\Omega\subset\mathbb R^3\to\mathbb S^2=\{y\in\R^3: |y|=1\}$, representing
mean orientations of molecule's optical axis, that minimizes the
Oseen-Frank bulk energy functional ${\calE}_{\rm{OF}}(d)=\int_\Omega \calW(d,\nabla d)\,dx$, where
\begin{eqnarray}\label{OF_density}
2\calW(d,\nabla d)&=&k_1({\rm{div}}\ d)^2+k_2(d\cdot{\rm{curl}}\ d)^2+k_3|d\times {\rm{curl}}\ d|^2\nonumber\\
&&+(k_2+k_4)[{\rm{tr}}(\nabla d)^2-({\rm{div}}\ d)^2],
\end{eqnarray}
where $k_1,k_2,k_3>0$ are splaying, twisting, and bending constants, $k_2\ge |k_4|$, and $2k_1\ge k_2+k_4$.
Utilizing the null-Lagrangian property of last term in $\calW(d,\nabla d)$, it can be shown that
under strong anchoring condition (or Dirichlet boundary condition) $d\big|_{\partial\Omega}
=d_0\in H^1(\Omega,\mathbb S^2)$, there
always exists a minimizer $d\in H^1(\Omega, \mathbb S^2)$ of the Oseen-Frank energy functional ${\calE}_{\rm{OF}}$. Furthermore,
such a $d$ solves the Euler-Lagrange equation:
\begin{equation}\label{OF_equation}
\frac{\delta\calW}{\delta d}\times d:=\big\{-{\rm{div}}(\frac{\partial\calW}{\partial \nabla d}(d,\nabla d))+\frac{\partial\calW}{\partial d}(d,\nabla d)\big\}\times d=0.
\end{equation}
The basic question for such a minimizer $d$ concerns both the regularity property and the defect structure. A fundamental result
by Hardt, Lin, and Kinderlehrer \cite{HLK} in 1986's (see also \cite{HKL1}) asserts
\begin{Theorem} \label{HLK1}
If $d\in H^1(\Omega, \mathbb S^2)$ is a minimizer of the Oseen-Frank energy ${\calE}_{\rm{OF}}$, then $d$ is analytic on $\Omega\ \setminus\ {\rm{sing}}(d)$
for some closed subset ${\rm{sing}}(d)\subset\Omega$, whose Hausdorff dimension is smaller than one.  If, in addition, $d_0\in C^{k,\alpha}(\overline\Omega,\mathbb S^2)$ for some $k\in\mathbb N_+$, then there exists a closed subset $\Sigma_1\subset\partial\Omega$, with $\calH^1(\Sigma_1)=0$,
such that $d\in C^{k,\alpha}(\overline\Omega\ \setminus\ ({\rm{sing}}(d)\cup\Sigma_1), \mathbb S^2)$.
\end{Theorem}

Concerning the interior singular set ${\rm{sing}}(d)$ in Theorem 1.1, an outstanding open question is
\begin{Question} {\it How large is the singular set ${\rm{sing}}(d)$ for minimizers of Oseen-Frank energy
$\calE_{\rm{OF}}$? or equivalently, what is the optimal estimate of size of $
{\rm{sing}}(d)$? }
\end{Question}

 Notice that when $k_1=k_2=k_3=1$ and $k_4=0$, the Oseen-Frank energy density function
$\calW(d,\nabla d)=\frac12|\nabla d|^2$ is  the Dirichlet energy density. Hence the problem reduces
to minimizing harmonic maps into $\mathbb S^2$ and (\ref{OF_equation}) becomes the equation of harmonic maps:
\begin{equation}\label{HM}
\frac{\delta\calW}{\delta d}:=\Delta d+|\nabla d|^2 d=0.
\end{equation}
Harmonic maps have been extensively studied in the past several decades, which are much better understood.
See \cite{SU1, SU2, SU3, BCL, LW1}. In particular, it is well-known that in dimension $n=3$, the singular set
of minimizing harmonic maps is at most a nonempty set of finitely many points (see \cite{SU1}). A typical example is the following
hedgehog singularity.
\begin{Example}  (\cite{BCL} \cite{Lin1}).
$\frac{x}{|x|}: \mathbb R^3\to\mathbb S^2$ is a minimizing harmonic map.
\end{Example}

Motivated by the isolated point singularity for minimizing harmonic maps, it is natural to ask

\begin{Question} {\it Is the set ${\rm{sing}}(d)$ for a minimizer of Oseen-Frank energy $\calE_{\rm{OF}}$
in Theorem 1.1 a set of finite points?}
\end{Question}

The main difficulty to investigate the size of singular set of minimizers $d$ to
the Oseen-Frank energy functional $\calE_{\rm{OF}}$ is that it is an open question whether an
energy monotonicity inequality, similar to that of harmonic maps \cite{SU1},
holds for $d$. In fact, the following seemly simple question remains open.
\begin{Question} {\it Assume $k_4=-k_2$ and $\displaystyle\max_{1\le i\le 3} |k_i-1|<<1$.
Suppose $d\in H^1(\Omega,\mathbb S^2) $ is a minimizer of the Oseen-Frank energy ${\calE}_{\rm{OF}}$.
Is the following monotonicity inequality true?
\begin{equation}\label{monotonicity}
r^{-1}\int_{B_r}|\nabla d|^2+\int_{B_R\setminus B_r} |x|^{-1}\big|\frac{\partial d}{\partial |x|}\big|^2
\le \big(1+\omega(R-r)\big) R^{-1}\int_{B_R}|\nabla d|^2
\end{equation}
for $0<r\le R$, with $\displaystyle\lim_{\tau\downarrow 0}\omega(\tau)=0$.}
\end{Question}
Notice that when $k_4=-k_2$ and $\displaystyle\max_{1\le i\le 3}|k_i-1|$ is sufficiently close to $0$, it has been
shown by \cite{HKL1} that ${\rm{sing}}(d)$ for any minimizer $d$ of $\calE_{\rm{OF}}$ has locally finite points.
The readers can consult with the survey article by Hardt-Kinderlehrer \cite{Hardt-Kinder} for some
earlier developments of Oseen-Frank theory. The book by Virga \cite{Virga} and the lecture note by Brezis \cite{Brezis} provide detailed exposition of Oseen-Frank theory.
There is an important survey article by Lin-Liu \cite{LL0} on the development of mathematical analysis of 
liquid cyrstals up to 2001. 
Partly motivated by these mathematical studies and physical experiments where line and surface defects
of nematic liquid crystals have been observed (see \cite{Kleman}), Ericksen \cite{Ericksen2} proposed the
so-called Ericksen model of liquid crystals with variable degree of orientation. It assumes that the bulk energy density depends on both the orientation vector $d$, with $|d|=1$, and orientational order $s\in [-\frac12, 1]$, which is given by
\begin{equation}\label{ericksen_density}
w(s, d, \nabla s, \nabla d)=w_0(s)+w_2(s, d,\nabla s, \nabla d),
\end{equation}
where $w_0(s)$ is roughly a W-shape potential function with $w_0(-\frac12)=w_0(1)=+\infty$.
The term $w_2$ is given by
$$2w_2(s,\nabla s, d,\nabla d)=2\calW(d,\nabla d)+K_5|\nabla s-(\nabla s\cdot d)d-\alpha (\nabla d)d|^2+K_6(\nabla s\cdot d-\beta {\rm{div}d})^2,$$
where $\calW(d,\nabla d)$ is given by (\ref{OF_density}), and $K_5>0, K_6>0, \alpha,\beta$ are coefficients
that may depend on $s$. We can view the Oseen-Frank model as a special case of Ericksen's model
by imposing the constraint on the orientational order $s=s^*$, with $w_0(s^*)=\min\{w_0(s): s\in [-\frac12,1]\}$.
The simplest form of Ericksen's bulk energy density is
\begin{equation}\label{Ericksen_density1}
w(s,d)=k|\nabla s|^2+s^2|\nabla d|^2+w_0(s).
\end{equation}
One can translate the problem of minimizers $(s,d)$ of Ericksen's energy functional (\ref{Ericksen_density1})
into the problem of minimizing harmonic maps $(s, u)$, with potentials $w_0$,
\begin{equation}\label{Ericksen_density2}
\min\int_\Omega \big((k-1)|\nabla s|^2+|\nabla u|^2 +w_0(s)\big)\,dx
\end{equation}
subject to the constraint $|s|=|u|$ and $(s,u)=(s_0, u_0)$ on $\partial\Omega$.
Then the problem (\ref{Ericksen_density2}) is essentially a minimizing harmonic map
into a circular cone $C_k=\{(s,u)\in\mathbb R\times\mathbb R^3:
|s|=\sqrt{k-1} |u|\}\subset\mathbb R^4$ for $k>1$ or $C_k=\{(s,u)\in\mathbb R\times\mathbb R^3:
|s|=\sqrt{1-k}|u|\}\subset \mathbb R^{3,1}$-the Minkowski space for $0<k<1$, Lin has proved
in \cite{LFH} (see also Ambrosio
\cite{Ambrosio1, Ambrosio2}, Ambrosio-Virga \cite{Ambrosio-Virga}, Hardt-Lin \cite{Hardt-Lin},
Lin-Poon \cite{Lin-Poon}, and Calderer-Golovaty-Lin-Liu \cite{CGLL} for further related works)
\begin{Theorem}\label{line_defect} For any $(s_0, u_0)\in H^\frac12(\partial\Omega, C_k)$, there exists a minimizer $(s,u)$ of (\ref{Ericksen_density1}) such that the following properties hold:
\begin{itemize}
\item[(1)] If $k>1$, then $(s,u)\in C^\alpha(\Omega)$ for some $0<\alpha<1$, and
$(s,u)$ is analytic away from $s^{-1}\{0\}$. If $s\not\equiv 0$, then $s^{-1}\{0\}$ has
Hausdorff dimension at most $1$.
\item[(2)] If $0<k<1$, then $(s,u)\in C^{0,1}(\Omega)$, and
$(s,u)$ is analytic away from $s^{-1}\{0\}$. If $s\not\equiv 0$,
then $s^{-1}\{0\}$ has Hausdorff dimension at most $2$.
\end{itemize}
\end{Theorem}

Since a liquid crystal molecule is indistinguishable between its head and tail, i.e. $d\approx -d$,  its odd order moments must vanish and the second order moment is believed to be most important to the energy. In general, one can
use a class of $3\times 3$ symmetric, traceless matrices $Q$ to describe the second order moments:
$$\calM=\big\{Q\in \R^{3\times 3}: \ Q=\int_{\mathbb S^2} (p\otimes p-\frac13\mathbb I_3)\,d\mu(p),
\ \mu\in \calP(\mathbb S^2)\big\},
$$
where $\calP(\mathbb S^2)$ denotes the space of probability measures on $\mathbb S^2$. It is clear that for any
order parameter tensor $Q\in\calM$, its three eigenvalues $\lambda_i(Q)\in [-\frac13, \frac23]$ ($1\le i\le 3$) and $\displaystyle\sum_{i=1}^3\lambda_i(Q)=1$. If $Q$ has two equal eigenvalues, it becomes a uniaxial nematic
or cholesteric phase so that it can be expressed as $Q=s\big(d\otimes d-\frac13 \mathbb I_3)$, with $|d|=1$ and $s\in [-\frac12, 1]$. If $Q\in\calM$ has three distinct eigenvalues, then the liquid material is biaxial.
The Landau-de Gennes bulk energy for nematic liquid crystal material was derived by de Gennes \cite{LDG1, LDG2}. A simplified form is given by
\begin{eqnarray}\label{DG_energy}
\calE_{\rm{DG}}(Q)&=&\int_\Omega \Big[\frac12\big(L_1|\nabla Q|^2+L_2({\rm{div}}Q)^2+L_3{\rm{tr}}(\nabla Q)^2+L_4\nabla Q\otimes\nabla Q:Q\big)
\nonumber\\
&&\quad+\big(\frac12A{\rm{tr}}(Q^2)-\frac13 B{\rm{tr}}(Q^3)+\frac14 C {\rm{tr}}(Q^2)^2\big)\Big].
\end{eqnarray}
It is well-known that if $Q$ is uniaxial, then the Landau-de Gennes energy can be essentially
reduced to  the Ossen-Frank energy
or the Ericksen energy. We would like to mention that it is not hard to show that there exists at least one minimizer $Q\in \calM$ of the Landau-de Gennes energy functional $\calE_{\rm{DG}}$, provided we
remove the constraint on the eigenvalues $-\frac13\le \lambda_i(Q)\le \frac23$, $1\le i\le 3$.
Moreover, the defect set of biaxial minimizers $Q\in\calM$ of the Landau-de Gennes energy
$\calE_{\rm{DG}}$ is contained in the set on which at least two of its eigenvalues are equal. The study of
such defect sets is extremely challenging and remains largely open.  Here we mention that there have been
several recent interesting works on the static case of Landau-de Gennes' model by J. Ball and his collaborators, and others.
See Ball-Zarnescu \cite{BZ1, BZ2}, Ball-Majumdar \cite{BM1}, Majumdar \cite{M}, Majumdar-Zarnescu \cite{MZ},
Henao-Majumdar \cite{HM}, Ngugen-Zarnescu \cite{NZ}, Bauman-Park-Phillips \cite{BPP},
and Evans-Kneuss-Tran
\cite{EKT}.

\section{Hydrodynamic Theory}

\subsection{Ericksen-Leslie system modeling nematic liquid crystal flows}
The Ericksen-Leslie system modeling the hydrodynamics of nematic liquid crystals, reducing
to the Oseen-Frank theory in the static case,  was proposed
by Ericksen and Leslie during the period between 1958 and 1968 (\cite{Ericksen} and \cite{Leslie}).
It is a macroscopic continuum description of the time evolution of the materials under the influence
of both the flow velocity field $u(x,t)$ and the macroscopic description of the microscopic orientation
configuration $d(x,t)$ of rod-like liquid crystals, i.e., $d(x,t)$ is a unit vector in $\R^3$. The full Ericksen-Leslie
system can be described as follows. The hydrodynamic equation takes the form
\begin{equation}\label{u-eqn}
\partial_t u+u\cdot\nabla u+\nabla P=\nabla\cdot\sigma,
\end{equation}
where the stress $\sigma$ is modeled by the phenomenological constitutive relation:
$$\sigma=\sigma^L(u,d)+\sigma^E(d).$$
Here $\sigma^L(u,d)$ is the viscous (Leslie) stress given by
\begin{eqnarray}\label{leslie}
\sigma^L(u,d)&=&\mu_1(d\otimes d:A)d\otimes d+\mu_2 d\otimes N+\mu_3N\otimes d
+\mu_4 A\nonumber\\
&&+\mu_5 d\otimes(A\cdot d)+\mu_6 (A\cdot d)\otimes d,
\end{eqnarray}
with six Leslie coefficients $\mu_1,\cdots,\mu_6$, and
$$A=\frac12(\nabla u+(\nabla u)^T, \ \omega=\frac12({\nabla u-(\nabla u)^T}),
\ N=\partial_t d+u\cdot\nabla d-\omega\cdot d$$
representing the velocity gradient tensor,  the vorticity field, and rigid rotation part of the changing rate of
the director by fluid respectively. Assuming the fluid is incompressible, we have
\begin{equation}\label{incompress}
\nabla\cdot u=0.
\end{equation}
While $\sigma^E(d)$ is the elastic (Ericksen) stress
$$\sigma^E(d)=-\frac{\partial\calW}{\partial(\nabla d)}\cdot (\nabla d)^T$$
where $\calW=\calW(d,\nabla d)$ is the Oseen-Frank energy density given by (\ref{OF_density}).  The dynamic equation for the director
field takes the form
\begin{equation}\label{d-eqn}
d\times\big(-\frac{\delta\calW}{\delta d}-\lambda_1 N -\lambda_2 A\cdot d\big)=0,
\end{equation}
where $\frac{\delta\calW}{\delta d}$ is the Euler-Lagrangian operator given by (\ref{OF_equation}).
We also have the compatibility condition
\begin{equation}\label{compat}
\lambda_1=\mu_2-\mu_3, \ \lambda_2=\mu_5-\mu_6,
\end{equation}
and Parodi's condition
\begin{equation}\label{Parodi}
\mu_2+\mu_3=\mu_6-\mu_5,
\end{equation}
which provides certain dynamical stability of the whole system 
(see \cite{parodi}).

\subsection{Energetic Variational Approach}  The
derivation of Ericksen-Leslie system in \cite{Ericksen} and \cite{Leslie}
is based on the conservation of mass, the incompressibility of
fluid, and conservation of both linear and angular momentums. Here we will sketch a modern approach
called energy variational approach, prompted by C. Liu and his collaborators (see \cite{SL}, \cite{WXL}, and \cite{LWX}),
which is based
on both the least action principle and the maximal dissipation principle.  In the context of hydrodynamics of nematic liquid crystals,
the basic variable is the flow map $x=x(X,t):\Omega_0^X\to\Omega_t^x$, with $X$ the Lagrangian (or material) coordinate and $x$
the Eulerian (or reference) coordinate. For a given velocity field $u(x,t)$, the flow map $x(X,t)$ solves:
$$x_t(x(X,t), t)=u(x(X,t),t); \ x(X,0)=X\in \Omega_0^X.$$
The deformation tensor ${\bf F}$ of the flow map $x(X,t)$ is
${\bf F}=\big(\frac{\partial x_i}{\partial X_j}\big)_{1\le i, j\le 3}$. By the chain rule, ${\bf F}$ satisfies the following
transport equation:
$$\partial_t{\bf F}+u\cdot\nabla {\bf F}=\nabla u {\bf F}.$$
The kinematic transport of the director field $d$ represents the molecules moving in the flow \cite{Larson, LWX}.
It can be expressed by
\begin{equation}\label{trans1}
d(x(X,t),t)={\bf E} d_0(X),
\end{equation}
where the deformation tensor ${\bf E}$ carries information of micro structures and configurations. It satisfies
the following transport equation:
\begin{equation}\label{trans2}
\partial_t{\bf E}+u\cdot\nabla{\bf E}=\omega{\bf E}+(2\alpha-1)A {\bf E},
\end{equation}
where $2\alpha-1=\frac{r^2-1}{r^2+1}\in [-1,1]$ ($r\in\R$) is related to the aspect ratio
of the ellipsoid shaped liquid crystal molecules, see Jeffrey \cite{Jeffrey}. Notice that if there is no internal damping, then the transport
equation of $d$ from (\ref{trans1}) and (\ref{trans2}) is given by
\begin{equation}\label{trans3}
\partial_t d+v\cdot\nabla d-\omega \cdot d-(2\alpha-1)A d=0.
\end{equation}

The total energy of the liquid crystal system is the sum of kinetic energy and internal elastic energy and is given by
$$\calE^{\rm{total}}=\calE^{\rm{kinetic}}+\calE^{\rm{int}},
\ \calE^{\rm{kinetic}}=\frac12\int |v|^2, \ \calE^{\rm{int}}={\calE}_{\rm{OF}}(d)=\int \calW(d,\nabla d).$$
The action functional $\calA$ of particle trajectories in terms of $x(X,t)$ is given by
$$\calA(x)=\int_0^T\big(\calE^{\rm{kinetic}}-\calE^{\rm{int}}\big)\,dt.$$
Since the fluid is assumed to be incompressible,  the flow map $x(X,t)$ is volume preserving, which
is equivalent to $\nabla\cdot u=0$.
The least action principle asserts that the action functional $\calA$ minimizes among all volume preserving
flow map $x(X,t)$, i.e., $\delta_x\calA=0$ subject to the constraint $\nabla\cdot u=0$. To simplify the derivation,
we consider one constant approximation of Oseen-Frank energy density function and
relax the condition $|d|=1$ by setting $\calW(d,\nabla d)=\frac12|\nabla d|^2+F(d)$ for $d:\Omega\to\R^3$.
Then direct calculations (see \cite{LWX} 7.1) yield that
\begin{equation}\label{LAP1}
0=\delta_x\calA=\frac{d}{d\epsilon}\big|_{\epsilon=0}\calA(x^\epsilon)
=\int_0^T\int_{\Omega}\big(\partial_t u+u\cdot\nabla u+\nabla\cdot(\nabla d\odot\nabla d)-\nabla\cdot
\widetilde\sigma\big) y\,dxdt,\footnote{$\nabla d\odot\nabla d=\big(\langle \frac{\partial d}{\partial x_i},
\frac{\partial d}{\partial x_j}\rangle\big)_{1\le i, j\le 3}$}
\end{equation}
where $y=\frac{d}{d\epsilon}\big|_{\epsilon=0} x^\epsilon$ satisfies $\nabla\cdot y=0$, and
$$\widetilde\sigma=-\frac12\big(1-\frac{\lambda_2}{\lambda_1}\big)(\Delta d-f(d))\otimes d+
\frac12\big(1+\frac{\lambda_2}{\lambda_1}\big)d\otimes (\Delta d-f(d)), \ f(d)=\nabla F(d).$$
Hence it follows from (\ref{LAP1}) that
\begin{equation}\label{LAP2}
\partial_t u+u\cdot\nabla u+\nabla P=-\nabla\cdot(\nabla d\odot\nabla d)+\nabla\cdot\widetilde\sigma,
\end{equation}
where the pressure $P$ serves as a Lagrangian multiplier for the incompressibility of fluid.

It follows from (\ref{trans2}) that the total transport equation of $d$, without internal microscopic damping, is
\begin{equation}\label{trans4}
\partial_t d+u\cdot\nabla d-\omega d+\frac{\lambda_2}{\lambda_1}Ad=0.
\end{equation}
If we take the internal microscopic damping into account, then we have
\begin{equation}\label{trans5}
\partial_t d+u\cdot\nabla d-\omega d+\frac{\lambda_2}{\lambda_1}Ad=\frac{1}{\lambda_1}
\frac{\delta \calE^{\rm{int}}}{\delta d}=-\frac{1}{\lambda_1}(\Delta d-f(d)).
\end{equation}
It is known that the dissipation functional $\calD$ to the system (\ref{SEL1}) is
\begin{eqnarray}\label{entropy1}
&&\calD=\int_\Omega \big[\mu_1 |d^TAd|^2+\frac12{\mu_4}|\nabla u|^2+(\mu_5+\mu_6)|Ad|^2\nonumber\\
&&\qquad\qquad+\lambda_1 |N|^2+(\lambda_2-\mu_2-\mu_3)\langle N, Ad\rangle\big].
\end{eqnarray}
By (\ref{trans5}) and Parodi's condition (\ref{Parodi}), $\calD$ can be rewritten as
\begin{eqnarray}\label{entropy2}
&&\calD=\int_\Omega \big[\mu_1 |d^TAd|^2+\frac12{\mu_4}|\nabla u|^2+(\mu_5+\mu_6+\frac{\lambda_2^2}{\lambda_1})|Ad|^2\nonumber\\
&&\qquad\qquad-\lambda_1 |N+\frac{\lambda_2}{\lambda_1}Ad|^2\big].
\end{eqnarray}

According to the maximal dissipation principle, one has that first order variation of $\calD$ with respect to
the rate function $u$ must vanish, i.e. $\frac{\delta\calD}{\delta u}=0$ subject to the constraint $\nabla\cdot u=0$. By direct calculations (see \cite{LWX} 7.2), we have that
$$
0=\frac{\delta\calD}{\delta u}=2\int_\Omega \langle u, \nabla\cdot(\nabla d\odot\nabla d)-\nabla\cdot\sigma^L(u,d)\rangle,
$$
where $\sigma^L(u,d)$ is the Leslie stress tensor given by (\ref{leslie}). Hence we arrive at
\begin{equation}\label{dissipation-force}
0=-\nabla P-\nabla\cdot(\nabla d\odot\nabla d)+\nabla\cdot\big(\sigma^{L}(u,d)\big).
\end{equation}
Combining the dissipative part derived from the maximal dissipation principle with the conservative
part derived from the least action principle, we arrive at the equation (\ref{u-eqn}). This, together with
(\ref{trans5}) and (\ref{incompress}), yields the Ginzburg-Landau approximated version of the Ericksen-Leslie system (\ref{u-eqn})-(\ref{incompress})-(\ref{d-eqn}).

\subsection{Brief descriptions on relationship between various models} There are basically three different kinds of theories to model nematic liquid crystals: Doi-Onsager theory, Landau-de Gennes theory, and Ericksen-Leslie theory.
The first is the molecular kinetic theory, and the latter two are the continuum theory.  Kuzzu-Doi \cite{KD} and
E-Zhang \cite{EZ} have formally derived the Ericksen-Leslie equation from the Doi-Onsager equation  by taking small Deborah number limit. Wang-Zhang-Zhang \cite{Wang-Zhang-Zhang1} have rigorously justified this formal derivation before the first singular time of the Ericksen-Leslie equation.  Wang-Zhang-Zhang \cite{Wang-Zhang-Zhang2} has rigorously derived Ericksen-Leslie equation from Beris-Edwards model in the
Landau-de Gennes theory. There are several  dynamic Q-tensor models describing hydrodynamic flows of
nematic liquid crystals, which are derived either by closure approximations (see \cite{FCL} \cite{FLS}
and \cite{Wang-Zhang-Zhang3}), or by variational methods such as Beris-Edwards' model \cite{BE} and Qian-Sheng's model \cite{QS}. See \cite{PZ1} and \cite{PZ2} for the well-posedness of some dynamic Q-tensor models.
Furthermore, a systematical approach to derive the continuum theory for nematic liquid crystals from the molecular kinetic theory in both static and dynamic cases was proposed in \cite{HLWZ} and \cite{Wang-Zhang-Zhang3}.
We would also like to mention the derivation of liquid crystal theory from the statistical view of points by
Seguin-Fried \cite{SF1, SF2}.

\subsection{Analytic issues on simplified Ericksen-Leslie system} When we consider  one constant approximation
for the Oseen-Frank energy, i.e., $\calW(d,\nabla d)=\frac12|\nabla d|^2$, the general Ericksen-Leslie system reduces
to
\begin{equation}\label{GEL2}
\begin{cases}
\partial_t u+u\cdot \nabla u+\nabla P=-\nabla\cdot(\nabla d\odot\nabla d)+\nabla\cdot(\sigma^L(u,d)),\\
\nabla\cdot u = 0,\\
\partial_t d+u\cdot\nabla d-\omega d+\frac{\lambda_2}{\lambda_1} Ad
=\frac{1}{|\lambda_1|}\left(\Delta d+|\nabla d|^2d\right)+\frac{\lambda_2}{\lambda_1}\left(d^TAd\right) d.
\end{cases}
\end{equation}
Direct calculations, using both (\ref{compat}) and (\ref{Parodi}), show that smooth solutions of (\ref{GEL2}), under suitable
boundary conditions, satisfy
the following energy law (see \cite{HLW}):
\begin{eqnarray}\label{global_energy_ineq1}
&&\frac{d}{dt}\int (|u|^2+|\nabla d|^2)
+\int \big[\mu_4|\nabla u|^2+\frac{2}{|\lambda_1|}|\Delta d+|\nabla d|^2d|^2\big]\nonumber\\
&&=-2\int \Big[\big(\mu_1-\frac{\lambda_2^2}{\lambda_1}\big)|A: d\otimes d|^2
+\big(\mu_5+\mu_6+\frac{\lambda^2_2}{\lambda_1}\big)|A\cdot d|^2\Big].
\end{eqnarray}
In particular, we have the following energy dissipation property.
\begin{Lemma}\label{energy_dissipation1} Assume both (\ref{compat}) and (\ref{Parodi}). If Leslie's coefficients
satisfy the algebraic condition:
\begin{equation}\label{algebra_cond}
\lambda_1<0, \ \mu_1-\frac{\lambda_2^2}{\lambda_1}\ge 0, \ \mu_4>0, \ \mu_5+\mu_6\ge -\frac{\lambda_2^2}{\lambda_1},
\end{equation}
then any smooth solution ($u,d$) of (\ref{GEL2}), under suitable
boundary conditions, satisfies the energy dissipation inequality:
\begin{equation}\label{energy_dissipation2}
\frac{d}{dt}\int (|u|^2+|\nabla d|^2)
+\int \big[\mu_4|\nabla u|^2+\frac{2}{|\lambda_1|}|\Delta d+|\nabla d|^2d|^2\big]\le 0.
\end{equation}
\end{Lemma}
A fundamental question related to the general Ericksen-Leslie system (\ref{GEL2}) is
\begin{Question} {\it For dimensions $n=2 \ {\rm{or}}\ 3$, smooth domains $\Omega\subset\mathbb R^n$ (or $\Omega=\mathbb R^n$), establish
\begin{itemize}
\item[i)] the existence of global Leray-Hopf type weak solutions of (\ref{GEL2}) under generic initial and boundary values $(u_0, d_0)\in {\bf H}\times H^1(\Omega,\mathbb S^2)$\footnote{Throughout this paper, ${\bf H}={\rm{Closure\ of}}
\ \big\{v\in C^\infty_0(\Omega,\mathbb R^n)\ | \ {\rm{div}} v=0\big\} \ {\rm{in}}\ L^2(\Omega,\mathbb R^n)$
and ${\bf J}={\rm{Closure\ of}}
\ \big\{v\in C^\infty_0(\Omega,\mathbb R^n)\ |\  {\rm{div}} v=0\big\} \ {\rm{in}}\ H^1_0(\Omega,\mathbb R^n)$.}.
\item[ii)] partial regularity and uniqueness properties for certain restricted classes of weak solutions of (\ref{GEL2}).
\item[iii)] an optimal global (or local) well-posedness of (\ref{GEL2}) for rough initial data ($u_0,d_0$), and the long time behavior
of global solutions of (\ref{GEL2}).
\end{itemize}}
\end{Question}

\noindent While both part i) and part ii) remain open for dimensions $n=3$ in general, the problem has been essentially solved in dimensions
$n=2$ recently. Furthermore, there has been some interesting progress towards the problem in dimensions $n=3$.

\smallskip
\noindent In order to rigorously analyze (\ref{GEL2}), Lin \cite{Lin2} first proposed a simplified version of (\ref{GEL2}) that preserves both the nonlinearity and the energy dissipation mechanism, i.e., $(u,d):\Omega\times (0,+\infty)\to \mathbb R^n\times \mathbb S^2$ solves
\begin{equation}\label{SEL1}
\begin{cases}
\partial_t u+u\cdot\nabla u+\nabla P=\mu \Delta u-\nabla\cdot(\nabla d\odot\nabla d),\\
\nabla\cdot u = 0, \\
\partial_t d+u\cdot\nabla d=\Delta d+|\nabla d|^2d.
\end{cases}
\end{equation}
It is readily seen that under suitable boundary conditions, (\ref{SEL1}) enjoys the following energy dissipation property:
\begin{equation}\label{energy_dissipation3}
\frac{d}{dt}\int (|u|^2+|\nabla d|^2)
=-2\int \big[\mu|\nabla u|^2+|\Delta d+|\nabla d|^2d|^2\big]\le 0.
\end{equation}
Because of the supercritical nonlinearity $\nabla\cdot(\nabla d\odot\nabla d)$ in (\ref{SEL1})$_1$,
Lin-Liu \cite{LL1, LL2, LL3} have studied the Ginzburg-Landau approximation (or orientation of variable degrees
in Ericksen's terminology \cite{Ericksen1}) of (\ref{SEL1}): for $\epsilon>0$, $(u,d):\Omega\times (0,+\infty)\to\mathbb R^n\times\mathbb R^3$
solves
\begin{equation}\label{SEL2}
\begin{cases}
\partial_t u+u\cdot\nabla u+\nabla P=\mu \Delta u-\nabla\cdot(\nabla d\odot\nabla d),\\
\nabla\cdot u = 0, \\
\partial_t d+u\cdot\nabla d=\Delta d+\frac{1}{\epsilon^2}\left(1-|d|^2\right)d.
\end{cases}
\end{equation}
Under the initial and boundary conditions:
\begin{equation}\label{IBC}
\begin{cases}
(u,d)\big|_{t=0}=(u_0,d_0), & x\in\Omega\\
(u,d)\big|_{\partial\Omega}=(0, d_0), & t>0,
\end{cases}
\end{equation}
we have the following basic energy law:
\begin{eqnarray}\label{energy_dissipation4}
&&\int_\Omega \big(|u|^2+|\nabla d|^2+\frac{1}{2\epsilon^2}(1-|d|^2)^2\big)(t)
+2\int_0^t\int_\Omega (\mu|\nabla u|^2+|\partial_t d+u\cdot\nabla d|^2)\nonumber\\
&&\le \int_\Omega \big(|u_0|^2+|\nabla d_0|^2+\frac{1}{2\epsilon^2}(1-|d_0|^2)^2\big).
\end{eqnarray}
Besides stability and long time asymptotic, the following existence and partial regularity were proven in \cite{LL1, LL2} (see also \cite{Liu-Walkington}).
\begin{Theorem} For any $\epsilon>0$ fixed, the following holds:
\begin{itemize}
\item[(a)]for $u_0\in {\bf H}$ and $d_0\in H^1(\Omega)\cap H^{\frac32}(\partial\Omega)$, there exists a global weak solution
$(u,d)$ of (\ref{SEL2}) and (\ref{IBC}) satisfying (\ref{energy_dissipation4}) and
$$\begin{cases} u\in L^2(0,T; {\bf J})\cap L^\infty(0,T; {\bf H}), \\
\ d\in L^2(0,T; H^1(\Omega))\cap L^\infty(0,T; H^2(\Omega)), \ \forall\ T\in (0, +\infty).
\end{cases}$$
\item[(b)] there exists a unique global classical solution $(u,d)$ to (\ref{SEL2}) and (\ref{IBC}) provided $(v_0,d_0)\in {\bf J}\times H^2(\Omega)$
and either $n=2$ or $n=3$ and $\mu\ge\mu(u_0,d_0).$
\item[(c)] if $(u,d)$ is a global suitable weak solution of (\ref{SEL2}), then $(u,d)\in C^\infty(\Omega\times (0,+\infty)\setminus\Sigma)$, with
$\calH^1(\Sigma)=0$.
\end{itemize}
\end{Theorem}
A long outstanding open problem pertaining to the Ginzburg-Landau approximation equation (\ref{SEL2}) and the original
nematic liquid crystal flow  equation (\ref{SEL1}) is
\begin{Question} {\it Whether weak limits $(u,d)$ of weak solutions $(u_\epsilon, d_\epsilon)$ of (\ref{SEL2}) and (\ref{IBC})
are weak solutions of (\ref{SEL1}) and (\ref{IBC}), as $\epsilon\rightarrow 0$.}
\end{Question}
It follows from (\ref{energy_dissipation4}) that there exists
$(u,d)\in \big(L^2(0,T; {\bf J})\cap L^\infty(0,T, {\bf H})\big)\times L^\infty (0,T; H^1(\Omega,\mathbb S^2))$
such that, up to a subsequence,
$$\begin{cases}
(u_\epsilon, d_\epsilon)\rightharpoonup (u,d) \ {\rm{in}}\ L^2(0,T; L^2(\Omega))\times L^2(0,T; H^1(\Omega)),\\
e_\epsilon(d_\epsilon)\,dxdt:=\big(|\nabla d_\epsilon|^2+\frac{1}{2\epsilon^2}(1-|d_\epsilon|^2)^2\big)\,dxdt
\rightharpoonup \mu:=|\nabla d|^2\,dxdt+\eta,\\
(\nabla d_\epsilon\odot\nabla d_\epsilon)\,dxdt\rightharpoonup (\nabla d\odot\nabla d)\,dxdt+\mathcal{M},
\end{cases}
$$
for some nonnegative Radon measure $\eta$ and a positive semi-definite $n\times n$ symmetric
matrix-valued function $\mathcal{M}$ with each entry being a Radon measure, as convergence of Radon measures
in $\Omega\times [0,T]$. Both $\nu$ and $\mathcal{M}$ are called {\it defect measures}, and ${\mathcal M}<<\nu$.
It is readily seen that
\begin{equation}\label{SEL3}
\begin{cases}
\partial_t u+u\cdot\nabla u+\nabla P=\mu \Delta u-\nabla\cdot(\nabla d\odot\nabla d+{\mathcal M}),\\
\nabla\cdot u = 0, \\
\partial_t d+u\cdot\nabla d=\Delta d+|\nabla d|^2d.
\end{cases}
\end{equation}
Thus we need to study the following difficult question.
\begin{Question} {\it Let $\eta$ and ${\mathcal M}$ be the defect measures.
\begin{itemize}
\item[1)] How large are the supports of $\nu$ and ${\mathcal M}$?
\item[2)] Under what conditions, $\nu$ or ${\mathcal M}$ vanishes?
\end{itemize}}
\end{Question}

Before describing a very recent progress towards Question 2.4 and 2.5, we present slightly earlier works on
the simplified Ericksen-Leslie system (\ref{SEL1}) in dimension $n=2$ by F. Lin, J. Lin, and C. Wang \cite{LLW1},
through a different approach that
is based on priori estimates under small energy condition. More precisely, we have established in \cite{LLW1}
\begin{Lemma} \label{epsilon-regularity} For $n=2$ and $0<r\le 1$, let $P_r=B_r\times [-r^2,0]$
be the parabolic ball of radius $r$.
There exists $\epsilon_0>0$ such that if $u\in L^\infty_tL^2_x\cap L^2_tH^1_x(P_r,\mathbb R^2)$,
$d\in L^2_tH^2_x(P_r,\mathbb S^2)$, and
$P\in L^2(P_r)$ solves (\ref{SEL1}), satisfying
\begin{eqnarray}\label{small_energy1}
&&\Phi(u, d, P, r):=\|u\|_{L^4(P_r)}+\|\nabla u\|_{L^2(P_r)}+\|\nabla d\|_{L^4(P_r)}\nonumber\\
&&\qquad\qquad\qquad\ +\|\nabla^2 d\|_{L^2(P_r)}
+\|P\|_{L^2(P_r)}\le\epsilon_0,
\end{eqnarray}
then $(u,d)\in C^\infty(P_{\frac{r}2}, \R^2\times\mathbb S^2)$, and
\begin{equation}\label{priori_estimate}
\big\|(u,\nabla d)\big\|_{C^l(P_\frac{r}2)}\le C(l)\epsilon_0 r^{-l}, \ \ l\ge 0.
\end{equation}
\end{Lemma}

Approximating initial and boundary data $(u_0, d_0)$ by smooth $(u_0^\epsilon, d_0^\epsilon)$,
and employing Lemma \ref{epsilon-regularity} to the short time smooth solutions $(u^\epsilon, d^\epsilon)$
of (\ref{SEL1}) under the initial and boundary value ($u_0^\epsilon, d_0^\epsilon$), we have established
the following result on existence and uniqueness of (\ref{SEL1}) in \cite{LLW1} and Lin-Wang \cite{LW2}
respectively.
\begin{Theorem}\label{existence} For any $u_0\in {\bf H}$ and $d_0\in H^1(\Omega,\mathbb S^2)$,
with $d_0\in C^{2,\alpha}(\partial\Omega,\mathbb S^2)$ for some $\alpha\in (0,1)$,
\\
{\rm{1)}} there exists a global weak solution $u\in L^\infty(0, +\infty; {\bf H})\cap L^2(0,+\infty; {\bf J})$ and
$d\in L^\infty(0,+\infty; H^1(\Omega,\mathbb S^2)$ of (\ref{SEL1}) and (\ref{IBC}) such that the
following properties hold:
\begin{itemize}
\item[1a)] there exists a nonnegative integer $L$ depending only on $u_0, d_0$ and $0<T_1<\cdots<T_L<+\infty$ such that
$$(u,d)\in C^\infty\big(\Omega\times [(0,+\infty)\setminus\{T_i\}_{i=1}^L]\big)
\cap C^{2,1}_\alpha\big(\overline\Omega\times[(0,+\infty)\setminus\{T_i\}_{i=1}^L]\big).$$
\item[1b)] at each time singularity $T_j$, $j=1,\cdots, L$, it holds
$$
\liminf_{t\uparrow T_j}\max_{x\in\overline\Omega}\int_{\Omega\cap B_r(x)}(|u|^2+|\nabla d|^2)(y,t)\,dy\ge 8\pi,
\ \forall r>0.$$
There exist $x_m^j\rightarrow x_0^j\in\Omega$ and $t_m^j\uparrow T_j$, $r_m^j\downarrow 0$, and
a nontrivial harmonic map $\omega_j:\R^2\to\mathbb S^2$ with finite energy such that
\begin{eqnarray*}
&&(u_m^j, d_m^j):=\big(r_m^j u(x_m^j+r_m^j x, t_m^j+(r_m^j)^2), d(x_m^j+r_m^j x, t_m^j+(r_m^j)^2)\big)\\
&&\rightarrow (0, \omega_j) \ {\rm{in}}\  C^2_{\rm{loc}}(\R^2\times [-\infty, 0]).
\end{eqnarray*}
\item[1c)] there exist $t_k\uparrow +\infty$ and a harmonic map $d_\infty\in C^\infty(\Omega,\mathbb S^2)\cap
C^{2,\alpha}_{d_0}(\overline\Omega,\mathbb S^2)$ such that
$u(t_k)\rightarrow 0$ in $H^1(\Omega)$, $d(t_k)\rightharpoonup d_\infty$ in $H^1(\Omega)$, and
there exist $\{x_i\}_{i=1}^l\subset\Omega$ and $\{m_i\}_{i=1}^l\subset \mathbb N$ such that
$$|\nabla d(t_k)|^2\,dx\rightharpoonup |\nabla d_\infty|^2\,dx+\sum_{i=1}^l 8\pi m_i \delta_{x_i}.$$
\item[1d)] if either $d_0^3\ge 0$ or $\int_\Omega (|u_0|^2+|\nabla d_0|^2)\le 8\pi$, then
$(u,d)\in C^\infty(\Omega\times (0,+\infty))\cap C^{2,1}_\alpha(\overline\Omega\times (0,+\infty))$.
Furthermore, there exists $t_k\uparrow +\infty$ and a harmonic map $d_\infty\in C^\infty(\Omega)\cap
C^{2,\alpha}_{d_0}(\overline\Omega,\mathbb S^2)$ such that $(u(t_k), d(t_k))\rightarrow (0,d_\infty)$
in $C^2(\overline\Omega)$.
\end{itemize}

\noindent {\rm{2)}} the global weak solution $(u,d)$ obtained in 1) is unique in the same class of weak solutions,
i.e. if $(v,e)\in \big(L^\infty(0,T_1^{-}; {\bf H})\cap L^2(0,T_1^{-}; {\bf J})\big)\times L^2(0,T_1^-; H^2(\Omega))$ is another
weak solution of (\ref{SEL1}) and (\ref{IBC}), then $(v,e)\equiv (u,d)$ in $\Omega\times [0,T_1)$.
\end{Theorem}

Some related works on the existence of global weak solutions of (\ref{SEL1}) in $\mathbb R^2$ have also been considered by
Hong \cite{Hong}, Hong-Xin \cite{HX},
Xu-Zhang \cite{XZ}, and Lei-Li-Zhang \cite{LLZ}.  There are a few interesting questions related to Theorem \ref{existence}, namely,
\begin{Question}{\it
1) Does the global weak solution $(u,d)$ obtained in Theorem \ref{existence} have at most finitely many
singularities?\\
2) Does there exist a smooth initial value $(u_0,d_0)$ such that the short time smooth solution $(u,d)$ to (\ref{SEL1}) and (\ref{IBC})
develop finite time singularity?\\
3) Does the solution $(u,d)$ obtained in Theorem \ref{existence} have a unique limit $(0, d_\infty)$ as $t\uparrow +\infty$?}
\end{Question}
It is not hard to check that the example of finite time singularity of harmonic map heat flow in dimension two by Chang-Ding-Ye \cite{CDY}
satisfies (\ref{SEL1}) with $u\equiv 0$. It is desirable to construct an example with a nontrivial velocity field. There has been
some partial result on the uniform limit at $t=\infty$ for (\ref{SEL1}) on $\mathbb S^2$, see \cite{WX}.

\begin{Remark} For the simplified Ericksen-Leslie system (\ref{SEL1}) in $\R^n$ for $n\ge 3$,
the uniqueness was proven for \\
i) the class of weak solutions $(u,d)\in C([0,T], L^n(\R^n,\R^n)\times W^{1,n}_{e_0}(\R^n,\mathbb S^2)$,
$e_0\in\mathbb S^2$, by Lin-Wang \cite{LW2}; and\\
ii) the class of weak solutions $(u,d) \in L^p(0,T; L^q(\R^n,\R^n))
\times L^p(0,T; W^{1,q}_{e_0}(\R^n, \mathbb S^2))$, with $p>2$ and $q>n$ satisfying Serrin's condition
$\frac{2}{p}+\frac{n}q=1$,   by Huang \cite{Huang}(it is interesting to ask whether the uniqueness holds  for the end point case
$(p,q)=(+\infty,n)$).
\end{Remark}

It turns out that the $\epsilon_0$-regularity lemma \ref{epsilon-regularity}  for the simplified Ericksen-Leslie
system (\ref{SEL1}) can also be proved by a blow-up type argument. Furthermore, with some delicate
analysis to establish smoothness of the limiting linear coupling system resulting from the blow-up argument,
such a blow up argument also works for the general Ericksen-Leslie system (\ref{GEL2}).
Indeed,  it was shown by Huang-Lin-Wang in \cite{HLW} that Lemma \ref{epsilon-regularity} also holds for (\ref{GEL2}).
As a consequence, we have extended in \cite{HLW} Theorem \ref{existence} to the general Ericksen-Leslie system
(\ref{GEL2}) in $\R^2$.
\begin{Theorem} \label{existence1} For $u_0\in {\bf H}$ and $d_0\in H^1_{e_0}(\R^2,\mathbb S^2)$, $e_0\in\mathbb S^2$,
assume the conditions (\ref{compat}), (\ref{Parodi}), and (\ref{algebra_cond}) hold. Then there exists a global
weak solution $(u,d):\R^2\times [0,+\infty)\to\R^2\times\mathbb S^2$ to (\ref{GEL2}) in $\R^2$, with
$u\in L^\infty(0,+\infty; {\bf H})\cap L^2(0,+\infty; {\bf J}), d\in L^\infty(0, +\infty; H^1_{e_0}(\R^2,\mathbb S^2))$, under the initial condition $(u,d)|_{t=0}=(u_0, d_0)$ in $\R^2$. Furthermore, ($u,d$)
enjoys the properties 1a), 1b), 1c), 1d) in Theorem \ref{existence} {\rm (}with $\Omega$ replaced by $\R^2${\rm )}.
\end{Theorem}

See also Wang-Wang \cite{WW} for some related work. It is an interesting question to ask whether
the uniqueness theorem for (\ref{SEL1}) by Lin-Wang \cite{LW2} also holds for
(\ref{GEL2}) in $\R^2$.

Now we outline a program to utilize the $\epsilon$-version of Ericksen-Leslie system (\ref{SEL2}) in the study of defect motion
of nematic liquid crystals in dimension $2$, see Lin \cite{Lin3} and \cite{Lin5}.   For a bounded smooth domain $\Omega\subset\R^2$,
let $(u_\epsilon, d_\epsilon):\Omega\times [0,+\infty)\to\R^2\times\R^2$ be solutions of (\ref{SEL2}),
and $d_\epsilon(x,t)=g(x): \partial\Omega\to\mathbb S^1$ is smooth with deg($g$)$=m\in \mathbb N$.
Assume the initial data $(u_0^\epsilon, d_0^\epsilon)$ satisfies
\begin{equation}\label{energy_bound}
\int_\Omega |u_0^\epsilon|^2\le {\bf A},\ E_\epsilon(d_0^\epsilon)
=\int_\Omega \big(|\nabla d_0^\epsilon|^2+\frac{1}{2\epsilon^2}(1-|d_0^\epsilon|^2)^2\big)\le 2\pi m\log \frac{1}{\epsilon}
+{\bf B}.
\end{equation}
Then the following facts hold:
\begin{itemize}
\item[(a)] $E_\epsilon(d_\epsilon(t))+\int_\Omega |u_\epsilon(t)|^2\le 2\pi m\log\frac{1}{\epsilon}+{\bf A}+{\bf B}$ for
$t\ge 0$.
\item[(b)] $E_\epsilon(d_\epsilon(t))\ge 2\pi m\log\frac{1}{\epsilon}-{\bf C}$,
 $\int_\Omega |u_\epsilon(t)|^2\le C({\bf A}, {\bf B}, {\bf C})$ for $t\ge 0$, and
$$\int_0^\infty\int_\Omega \big(\mu|\nabla u_\epsilon|^2+|\Delta d_\epsilon+\frac{1}{\epsilon^2}(1-|d_\epsilon|^2)d_\epsilon|^2\big)\le C({\bf A}, {\bf B}, {\bf C}).$$
\item[(c)] $\big\|\nabla d_\epsilon(t)\big\|_{L^p(\Omega)}\le C(p, \Omega, g, {\bf A}, {\bf B}, {\bf C})$ for $1\le p<2$.
\item[(d)] $u_\epsilon(t)\rightarrow u(t)$ in $L^2(\Omega)$ and weakly in $H^1(\Omega)$, and
$\exists$ $h(t)\in H^1(\Omega)$, with $\|h(t)\|_{H^1(\Omega)}\le C$,
and $\{a_j(t)\}_{j=1}^m\subset C^{0,1}([0, T), \Omega)$ such that, up to a subsequence,
$$d_\epsilon(t)\rightarrow \displaystyle\Pi_{j=1}^m \frac{x-a_j(t)}{|x-a_j(t)|} e^{ih(x,t)}=e^{i\Theta_a+h}$$
 in $L^2(\Omega)\cap H^1_{\rm{loc}}(\overline\Omega\ \setminus\ \{a_j(t)\}_{j=1}^m).$
\item[(e)] $\nabla d_\epsilon\odot\nabla d_\epsilon\rightharpoonup (\nabla\Theta_a+\nabla h)\odot (\nabla\Theta_a+\nabla h)+\eta$
in $\mathcal {D}'(\Omega)$ for some positive semi-definite matrix valued function $\eta$ with each entry being Radon measures in $\Omega$. Moreover, ${\rm{supp}}(\eta(t))\subset\{a_j(t)\}_{j=1}^m$.
\end{itemize}
Then we have
\begin{equation}\label{defect_motion}
\frac{d a_j}{dt}=u(a_j(t), t); \ a_j(0)=a_{0, j}, \ j=1, \cdots, m,
\end{equation}
and $(u, h, \Theta_a)$ satisfies
\begin{equation}\label{motion_eqn}
\begin{cases}
h_t+u\cdot \nabla h=\Delta h+(u(a(t), t)-u)\cdot\nabla\Theta_a,\\
u_t+u\cdot\nabla u+\nabla P=\mu\Delta u-\nabla h\cdot\Delta h -\nabla\Theta_a(\Delta h-\Delta h(a(t), t)),\\
\nabla\cdot u =0.
\end{cases}
\end{equation}
An important question is to show that the system (\ref{defect_motion}) and (\ref{motion_eqn}) admits a smooth solution
$(a(t), u, h, P)$. Once this has been established, the consistency of energy laws will yield the defect measure
$\eta\equiv 0$.

\medskip
Next we turn to discuss some recent progress towards the simplified Ericksen-Leslie system (\ref{SEL1}) in dimensions $n\ge 3$.
For (\ref{SEL1}), it is standard (see \cite{CS} \cite{Sh}) that for any initial data $(u_0,d_0)\in H^3(\mathbb R^n,\mathbb R^n)\times H^4_{e_0}(\mathbb R^n,\mathbb S^2)$, $e_0\in\mathbb S^2$, there exist $T_*>0$ and a unique strong solution $(u,d):\R^n\times [0,T_*)\to\R^n\times\mathbb S^2$ of (\ref{SEL1}) and (\ref{IBC}), see \cite{Hu-Wang} \cite{Li-Wang1} \cite{Li-Wang2}
for related works on the existence of global strong solutions to (\ref{SEL1}) for small initial data in suitable
function spaces. Furthermore, a criterion on possible breakdown for
local strong solutions of (\ref{SEL1}), analogous to BKM criterion for the Navier-Stokes equation \cite{BKM},  was obtained
by Huang-Wang \cite{Huang-Wang}: if $0<T_*<+\infty$ is the maximum time, then
\begin{equation}\label{BKM1}
\int_0^{T_*}\big(\|\nabla\times u\|_{L^\infty(\R^n)}+\|\nabla d\|_{L^\infty(\R^n)}^2\big)\,dt=+\infty.
\end{equation}
See Wang-Zhang-Zhang \cite{WZZ1} and Hong-Li-Xin \cite{HLX}
for related works on the general Ericksen-Leslie system (\ref{GEL2}).

For (\ref{SEL1}), inspired by Koch-Tataru \cite{KT} the following well-posedness result
has been obtained by Wang \cite{Wang1}.

\begin{Theorem}\label{WP1} There exists $\epsilon_0>0$ such that if $(u_0, d_0):\mathbb R^n\to\R^n\times\mathbb S^2$,
with $\nabla\cdot u_0=0$, satisfies
$$\big[u_0\big]_{{\rm{BMO}}^{-1}(\R^n)}+\big[d_0\big]_{{\rm{BMO}}(\R^n)}\le\epsilon_0,$$
then there exists a unique global smooth solution $(u,d)\in C^\infty(\R^n\times (0,+\infty),\R^n\times\mathbb S^2)$
of (\ref{SEL1}), under the initial condition $(u_0,d_0)$, which enjoys the decay estimate
\begin{equation}\label{decay_estimate}
t^\frac12\big(\|u(t)\|_{L^\infty(\R^n)}+\|\nabla d(t)\|_{L^\infty(\R^n)}\big)\lesssim \epsilon_0, \ t>0.
\end{equation}
\end{Theorem}

Higher order decay estimates for the solution given by Theorem \ref{WP1} have been established by \cite{JLin}
and \cite{Du}.  Local well-posedness of (\ref{SEL1}) in $\R^3$ for initial data $(u_0,d_0)$, with $(u_0,\nabla d_0)\in L^3_{\rm{uloc}}(\R^3)$ (the uniformly locally $L^3$-space in $\R^3$) having small norm
$\|(u_0,\nabla d_0)\|_{L^3_{\rm{uloc}}(\R^3)}$, has been
established by Hineman-Wang \cite{HW}. We would like to mention some interesting works on both new modeling and analysis of the hydrodynamics of non-isothermal nematic liquid crystals
by Feireisl-Rocca-Schimperna \cite{FRS}, Feireisl-Fr\'emond-Rocca-Schimperna \cite{FFRS1}, and Li-Xin
\cite{Li-Xin}.

Now we describe a very recent work by Lin and Wang \cite{LW3} on the existence of global
weak solutions of (\ref{SEL1}) in dimensions three by performing blow up analysis of the
Ginzburg-Landau approximation system (\ref{SEL2}). More precisely, we have proved

\begin{Theorem}\label{3d_existence} For $\Omega\subset\mathbb R^3$ either a bounded domain or the entire $\R^3$,
assume $u_0\in {\bf H}$ and $d_0\in H^1(\Omega,\mathbb S^2)$ satisfies $d_0(\Omega)\subset\mathbb S^2_+$, the upper half sphere. Then there exists a global weak solution $(u,d):\Omega\times [0,+\infty)\to\mathbb R^3\times \mathbb S^2$ to
the initial and boundary value problem of (\ref{SEL1}) and (\ref{IBC}) such that\\
(i) $u\in L^\infty(0, +\infty;{\bf H})\cap L^2(0,+\infty; {\bf J})$.\\
(ii) $d\in L^\infty(0,+\infty; H^1(\Omega,\mathbb S^2))$ and $d(t)(\Omega)\subset\mathbb S^2_+$ for
$L^1$ a.e. $t\in [0,+\infty)$.\\
(iii) $(u,d)$ satisfies the global energy inequality: for $L^1$ a.e. $t\in [0,+\infty)$,
\begin{equation}\label{global_energy_ineq}
\int_\Omega(|u|^2+|\nabla d|^2)(t)+2\int_0^t\int_\Omega \big(\mu|\nabla u|^2+|\Delta d+|\nabla d|^2 d|^2\big)
\le \int_\Omega(|u_0|^2+|\nabla d_0|^2).
\end{equation}
\end{Theorem}

The proof of Theorem \ref{3d_existence}  in \cite{LW3} is very delicate, which is based on suitable extensions of
the blow-up analysis scheme that has been developed in the context of harmonic maps  by Lin \cite{Lin4}
and harmonic map heat flows by Lin-Wang \cite{LW4, LW5, LW6}.  A crucial ingredient is the following compactness
result.
\begin{Theorem}\label{compactness}
For any $0<a\le 2$, $L_1>0$, and $L_2>0$, set $\mathbf{X}({L_1, L_2, a;\Omega})$ consisting of maps
$d_\epsilon\in H^1(\Omega,\mathbb R^3)$, $0<\epsilon\le 1$, that are solutions of
\begin{equation}\label{AGL5}
\Delta d_\epsilon+\frac{1}{\epsilon^2}(1-|d_\epsilon|^2) d_\epsilon=\tau_\epsilon\ \ {\rm{in}}\ \ \Omega
\end{equation}
such that the following properties hold:
\begin{itemize}
\item[(i)] $|d_\epsilon|\le 1$ and $d_\epsilon^3\ge -1+a$ for a.e. $x\in\Omega$.
\item[(ii)] $\displaystyle E_\epsilon(d_\epsilon)=\int_\Omega e_\epsilon(d_\epsilon)\,dx\le L_1$.
\item[(iii)] $\displaystyle\big\|\tau_\epsilon\big\|_{L^2(\Omega)}\le L_2$.
\end{itemize}
Then $\mathbf{X}(L_1, L_2,a; \Omega)$ is precompact in
$H^1_{\rm{loc}}(\Omega,\mathbb R^3)$, i.e., if $\{d_\epsilon\}\subset {\bf X}(L_1,L_2, a; \Omega)$,
then there exists $d_0\in H^1(\Omega,\mathbb S^2)$ such that, up to a subsequence, $d_\epsilon\rightarrow d_0$
in $H^1_{\rm{loc}}(\Omega, \mathbb R^3)$ as $\epsilon\rightarrow 0$.
\end{Theorem}

\medskip
\noindent{\bf \underline{Sketch of proof of Theorem \ref{compactness}}}.
Assume $d_\epsilon\rightharpoonup d_0$ in $H^1(\Omega)$ and
$$e_\epsilon(d_\epsilon)\,dx\rightharpoonup\mu:=\frac12|\nabla d_0|^2\,dx+\nu$$
as convergence of Radon measures  for some nonnegative Radon measure $\nu$, called defect measure.
We claim that $\nu\equiv 0$. This will be done in several steps.

\smallskip
\noindent{\it Step1} (almost monotonicity). $d_\epsilon\in {\bf X}(L_1,L_2, a;\Omega)$ satisfies
\begin{equation}\label{almost_mono1}
\Phi_\epsilon(R)\ge \Phi_\epsilon(r)+\int_{B_R\setminus B_r}|x|^{-1}\big|\frac{\partial d_\epsilon}{\partial |x|}\big|^2, \ \forall\ 0< r\le R,
\end{equation}
where
$$\Phi_\epsilon(r):=\frac1{r}\int_{B_r}\big(e_\epsilon(d_\epsilon)-\langle x\cdot\nabla d_\epsilon, f_\epsilon\rangle\big) +\frac12 \int_{B_r}|x||f_\epsilon|^2.$$

\medskip
\noindent{\it Step2} ($\delta_0$-strong convergence). $\exists \ \delta_0>0,$ $\alpha_0\in (0,1),$ and $C_0>0$ such that
if $d_\epsilon\in {\bf X}(L_1, L_2, a; \Omega)$ satisfies
$$\Phi_\epsilon(r_0)\le \delta_0,$$
then $|d_\epsilon|\ge \frac12$ in $B_{\frac{r_0}2}$, and
$$\Big[\frac{d_\epsilon}{|d_\epsilon|}\Big]_{C^{\alpha_0}(B_{\frac{r_0}2})}\le C_0.$$
In particular, $d_\epsilon\rightarrow d_0\ {\rm{in}}\ H^1(B_{\frac{r_0}2})$ and $\nu=0 \ {\rm{in}}\ B_{\frac{r_0}2}.$

\medskip
\noindent{\it Step3} (almost monotonicity for $\mu$)
\begin{equation}\label{almost_mono2}
\Theta^1(\mu, r):=\frac{1}{r}\mu(B_r)\le \Theta^1(\mu, R)+C_0(R-r), \ \forall \ 0<r\le R.
\end{equation}
Thus $\displaystyle\Theta^1(\mu)=\lim_{r\downarrow 0}\Theta^1(\mu, r)$ exists,
and is upper semicontinuous.

\medskip
\noindent{\it Step4} (concentration set). $d_\epsilon\rightarrow d_0$ in $H^1_{\rm{loc}}(\Omega\ \setminus\ \Sigma)$,
where
$$\displaystyle\Sigma:=\big\{x\in\Omega: \Theta^1(\mu)\ge \delta_0\big\}$$
is a $1$-rectifiable, closed subset, with $H^1(\Sigma)<+\infty$. Moreover,
$${\rm{supp}}(\nu)\subset\Sigma\ \ {\rm{and}}\ \delta_0\le \Theta^1(\nu, x)=\Theta^1(\mu,x)\le C_0 \   H^1{\rm{a.e.}}\ x\in\Sigma.$$

\medskip
\noindent {\it Step5} (stratification and blow-up). If $\nu\not\equiv 0$, then $H^1(\Sigma)>0$.
We can pick up a generic point $x_0\in\Sigma$ such that
$$\lim_{r_i\downarrow 0} r_i^{-1}\int_{B_{r_i}(x_0)}|\nabla d_0|^2=0,
\ \Theta^1(\nu, \cdot) \ {\rm{is}}\ H^1 {\rm{approximately \ continuous \ at }} \ x_0.$$
Define blow-up sequences $d_i(x)=d_{\epsilon_i}(x_0+r_i x)$ and $\nu_i(A)=r_i^{-1}\nu(x_0+r_iA)$ for $A\subset\R^3$.
Then
$$\nu_i\rightharpoonup \nu_0:=\theta_0 H^1{\rm{L}} Y , \ e_{\epsilon_i}(d_i)\,dx\rightharpoonup\nu_0
{\rm{(}}=\theta_0 H^1{\rm{L}} Y{\rm{)}}$$
as convergence of Radon measures, for some $\theta_0>0$ and $Y=\{(0,0,x_3): x_3\in\mathbb R\}$.
Applying (\ref{almost_mono1}), one has
$$\int_{B_2}\big|\frac{\partial d_{i}}{\partial x_3}\big|^2\rightarrow 0.$$
Now we can perform another round of blow up process to extract a nontrivial smooth harmonic map
$\omega:\mathbb R^2\to\mathbb S^{2}_{-1+\delta}$ that has finite energy (see \cite{LW3} for the
detail).
This is impossible, since any such a $\omega$ has non-zero topological degree.
\qed

\subsection{Compressible flow of nematic liquid crystals} When the fluid is assumed to compressible, the Ericksen-Leslie
system becomes more complicate and there seems very few analytic works available yet. We would like to mention that there have been both modeling study (see \cite{Morro}) and numerical study (see \cite{ZV}) on the hydrodynamics of compressible nematic liquid crystals under the influence of temperature gradient or electromagnetic forces. Here we briefly
describe some recent studies on a simplified version of compressible Ericksen-Leslie system in $\Omega\subset\R^n$,
which is given by
\begin{equation}\label{CNLF}
\begin{cases}
\partial_t \rho+\nabla \cdot (\rho u)=0,\\
\partial_t(\rho u)+\nabla\cdot(\rho u\otimes u)+\nabla (P(\rho))=\calL u-\nabla\cdot(\nabla d\odot\nabla d-\frac12|\nabla d|^2 \mathbb I_n),\\
\partial_t d+u\cdot\nabla d=\Delta d+|\nabla d|^2d,\\
\big(\rho, u, d\big)\big|_{t=0}=\big(\rho_0, u_0, d_0\big), \ \big(u,d\big)\big|_{\partial\Omega}=\big(u_0,d_0\big),
\end{cases}
\end{equation}
where $\rho:\Omega\subset\R^n\to\mathbb R_+$ is the fluid density and
$P(\rho):\Omega\to\mathbb R_+$ denotes the pressure of the fluid,
and $\calL$ denotes the Lam$\acute{\rm e}$ operator
$$
\calL u=\mu\Delta u+(\mu+\lambda)\nabla(\nabla\cdot u),
$$
where $\mu$ and $\lambda$ are the shear viscosity and the bulk viscosity coefficients
of the fluid, satisfying the following physical condition:
$$
\mu>0, \ 3\lambda+2\mu\geq 0.
$$
The system (\ref{CNLF}) is a strong coupling between compressible Navier-Stokes equation and the transported
harmonic map heat flow to $\mathbb S^2$. Since it is a challenging question to establish the existence of global
weak solutions to compressible Navier-Stokes equation itself (see \cite{Feireisl}), it turns out to be more difficult
to show the existence of global weak solutions to (\ref{CNLF}) than the incompressible nematic liquid crystal
flow (\ref{SEL1}). See \cite{DHXWZ} for a rigorous proof of convergence of (\ref{CNLF}) to (\ref{SEL1}) under
suitable conditions. The local existence of strong solutions to (\ref{CNLF}) can be established under suitable regularity
and compatibility conditions on initial data. For example, we have proven in \cite{HWW1, HWW2}
\begin{Theorem}\label{strong-solution} For $n=3$, assume $P\in C^{0,1}(\R_+)$,  $0\le\rho_0\in W^{1,q}\cap H^1\cap L^1$
for some $3<q<6$, $u_0\in H^1\cap H^2$, and $d_0\in H^2(\Omega, \mathbb S^2)$ satisfies
$$\calL u_0-\nabla(P(\rho_0))-\Delta d_0\cdot\nabla d_0=\rho_0^\frac12 g, \ g\in L^2.$$
Then there exist $T_0>0$ and a unique strong solution $(\rho, u, d)$ to (\ref{CNLF}) in $\Omega\times [0,T_0)$.
Furthermore, if $T_0<+\infty$ is the maximal time interval, then
\begin{equation}\label{BKM2}
\begin{cases}
\|\nabla u+(\nabla u)^t\|_{L^1(0, T_0; L^\infty)}+\|\nabla d\|_{L^2(0,T_0; L^\infty)}=+\infty,\\
{\rm{or}}\ \|\rho\|_{L^\infty(0,T_0; L^\infty)}+\|\nabla d\|_{L^3(0,T_0; L^\infty)}=+\infty
\ ( {\rm{if}}\ 7\mu>9\lambda).
\end{cases}
\end{equation}
\end{Theorem}
The local strong solution established in Theorem \ref{strong-solution} has been shown to be global for small
initial-boundary data $(\rho_0, u_0, d_0)$ by \cite{JSW} and \cite{LZZ}.

The global existence of  weak solutions to (\ref{CNLF}) is more challenging than (\ref{SEL1}) in dimensions $n\ge 2$,
see \cite{DLWW, DWW} for $n=1$. When relaxing the condition $|d|=1$ in (\ref{CNLF}) to $d\in\R^3$ by the
Ginzburg-Landau approximation, similar to (\ref{SEL2}),  the global existence of weak solutions
has been established by \cite{Wang-Yu} and \cite{LQ} in dimensions $n=3$. Finally, we would like to mention
that by employing the precompactness theorem \ref{compactness}, the following global existence of weak solutions
to (\ref{CNLF}) has been established by \cite{LLW}  in dimensions $n=2, 3$ (see \cite{JSW} a different proof for $n=2$).
\begin{Theorem}\label{weak-solution} Assume $P(\rho)=\rho^\gamma$ for $\gamma>1$,
$0\le \rho_0\in L^\gamma$, $m_0\in L^{\frac{2\gamma}{\gamma+1}}$ and
$\frac{|m_0|^2}{\rho_0}\chi_{\{\rho_0>0\}}\in L^1$, and $d_0\in H^1(\Omega, \mathbb S^2)$ with $d_0^3\ge 0$.
Then there exists a global renormalized weak solution $(\rho, u, d)$ to (\ref{CNLF}).
\end{Theorem}

\bigskip

\noindent{\bf Acknowledgements}. The first author is partially supported by NSF grants DMS1065964 and DMS1159313. The second author is partially supported by NSF grant DMS1001115 and DMS1265574, and
NSFC grant 11128102 and a Simons Fellowship in Mathematics. The second author wishes to thank his close friend
C. Liu for many thoughtful discussions on the subject.


\end{document}